%% file: mgr_multiphase.tex
\documentclass[article]{siamart}

\usepackage{amsmath}
\usepackage{amssymb}
\usepackage{float}
\usepackage{cleveref}
\usepackage{blkarray}
\usepackage{bm}
\usepackage{color}
\usepackage{tabularx}
\usepackage{caption}
\usepackage{multirow}
\usepackage{graphicx}
\usepackage[margin=1in]{geometry}
\usepackage[ruled, vlined]{algorithm2e}

\newcommand*{\mathcolor}{}
\def\mathcolor#1#{\mathcoloraux{#1}}
\newcommand*{\mathcoloraux}[3]{%
  \protect\leavevmode
  \begingroup
    \color#1{#2}#3%
  \endgroup
}

\input{ex_shared}

\ifpdf
\hypersetup{
  pdftitle={\TheTitle},
  pdfauthor={\TheAuthors}
}
\fi

\begin{document}
\maketitle

\begin{abstract}
Multiphase flow is a critical process in a wide range of applications, including oil and gas recovery, carbon sequestration, and contaminant remediation. Numerical simulation of multiphase flow requires solving of a large, sparse linear system resulting from the discretization of the partial differential equations modeling the flow. In the case of multiphase multicomponent flow with miscible effect, this is a very challenging task. The problem becomes even more difficult if phase transitions are taken into account. A new approach to handle phase transitions is to formulate the system as a nonlinear complementarity problem (NCP). Unlike in the primary variable switching technique, the set of primary variables in this approach is fixed even when there is phase transition. Not only does this improve the robustness of the nonlinear solver, it opens up the possibility to use multigrid methods to solve the resulting linear system. The disadvantage of the complementarity approach, however, is that when a phase disappears, the linear system has the structure of a saddle point problem and becomes indefinite, and current algebraic multigrid (AMG) algorithms cannot be applied directly. In this study, we explore the effectiveness of a new multilevel strategy, based on the multigrid reduction technique, to deal with problems of this type. We demonstrate the effectiveness of the method through numerical results for the case of two-phase, two-component flow with phase appearance/disappearance. We also show that the strategy is efficient and scales optimally with problem size.
\end{abstract}


\section{Introduction}
Modeling multiphase flow in porous media is a challenging task given the complex physics involved. The flow is described by a set of nonlinear and strongly coupled partial differential equations (PDEs) with algebraic constraints. These equations are usually solved with fully implicit method which requires solving large-scale, non-symmetric, and ill-conditioned linear systems.  The problem becomes even more difficult if we take into account phase transitions. When phase transitions occur, the PDEs can become degenerate and that makes the resulting linear system indefinite. Krylov subspace methods, such as the generalized residual method (GMRES) \cite{Saad86}, can be applied to solve these systems. However, these methods by themselves generally converge slowly and they must be appropriately preconditioned to accelerate convergence. The Incomplete LU (ILU) factorization is a popular approach as a preconditoner, due to its simplicity and generality. However, as simulations cover larger and larger domains and are deployed over high performance parallel architectures, there is an apparent need for robust solvers that scale, and the use of standard (single-level) ILU methods becomes less favorable.

Previously, most of the research has focused on finding new formulations which can deal with phase transitions. Some approaches include primary variable switching (PVS) \cite{Forsyth95,WuForsyth01}, negative saturation \cite{Abadpour09}, and finding a set of persistent primary variables \cite{Bourgeat09,Marchand12,Neumann13}. Recently, a new approach has been developed for handling the phase transitions by formulating the system of equations as a nonlinear complementarity problem (NCP) \cite{BenGharbia14,Lauser11,Marchand14}. Unlike the PVS approach, the advantage of the NCP approach is that the set of primary variables is consistent through out the simulation, and no primary variable switching is needed. Not only is this approach more robust and efficient, it also presents an opportunity to use scalable linear solvers such as algebraic multigrid (AMG). 

In this work, we focus on developing a new AMG preconditioner based on multigrid reduction (MGR) for GMRES to solve the linear system resulting from the discretization of the continuous problem. MGR technique has been around for many years \cite{Ries79,Ries83}. It can be considered as a generalization of the multi-stage preconditioner in a standarfd multigrid framework. A closed form of the error propagator can be derived for the MGR approach, and this enables us to study the effect of different multigrid components on the convergence of the linear solver. In addition, the MGR framework has been shown to be an efficient preconditioner for different types of PDEs, such as the reservoir simulation, and it has also been applied with varying degree of success to the time dimension \cite{Falgout14}. 

\textcolor{black}{We consider a two-phase, two-component system with phase transitions as our model problem. We describe this model in details in section \ref{sec:statement}. Classical approach to simulate two-phase, two-component is well-posed if two primary variables are chosen in advance. For example, one can choose one phase pressure and one phase saturation, or one phase pressure and one component concentration. This set of variables remains fixed in the case of phase appearance/disappearance if we know in advance what phase will appear and disappear during the simulation. Then a constrained pressure residual (CPR) preconditioning approach \cite{Dawson97,Qiao2017,Wallis85} can be employed to obtain a semi-elliptic pressure equation. The pressure equation is solved with a multi-level method, such as AMG or multiscale \cite{Cusini2015}, and followed by a relaxation step with ILU for the global linear system. Although this approach has been shown to be very effective for some real world examples \cite{Cusini2015}, it could be less robust in cases with strong capillarity effect \cite{Bui17}. Due to the fact that the CPR approach uses ILU in the smoothing step, it may not scale as well as a block factorization approach \cite{Bui17}. The goal of this paper is to develop a new multigrid algorithm that is both robust, efficient, and also general to accommodate various formulations of compositional multiphase flow. In particular, we show that, under appropriate assumptions, our multigrid reduction method is equivalent to the CPR-AMG and block factorization approaches.}

\textcolor{black}{The rest of this paper is organized as follows: In section \ref{sec:algorithm}, we describe the fully implicit discretization. We briefly review the MGR framework in section \ref{sec:MGR} and explain our new MGR algorithm in section \ref{sec:MGR-apply}. In section \ref{sec:result}, several numerical tests are studied for the robustness and scalability of the new algorithm. Some conclusion remarks as well as future work are presented in section \ref{sec:conclusion}.}

\section{Problem Statement}
\label{sec:statement}
\subsection{Governing Equations}
%\subsection{Simplified 2-phase, 2-component model}
We consider a simplified two-phase two-component model with phase transitions, similar to that presented in \cite{BourgeatBenchmark13}. This model provides a simple example that demonstrates the capability of the nonlinear complementarity constraint approach to handle phase appearance and disappearance. The flow consists of gas and liquid phases, and the components are hydrogen and water. We make the following simplifications: (1) water does not vaporize so the gas phase contains only hydrogen, and (2) the amount of hydrogen dissolved into the liquid phase is small. For a complete set of assumptions, we refer to \cite{Bourgeat09}. For the two components, the mass conservation equations read
\begin{align}
%&\phi \dfrac{\partial (\rho _l^w S_l)}{\partial t} - \nabla \cdot (\rho _l^w K\lambda _l\nabla P_l + j_l^h) = 0, \\
&\phi \dfrac{\partial (\rho _l^w S_l)}{\partial t} + \nabla \cdot (\rho _l^w \vec{q}_l - j_l^h) = 0 \label{water_mass_conservation},\\
%&\phi \dfrac{\partial (\rho _l^h S_l + \rho _g^h S_g)}{\partial t} - \nabla \cdot (\rho _l^h K\lambda _l\nabla P_l + \rho _g^h K\lambda _g\nabla P_g - j_l^h) = 0,
&\phi \dfrac{\partial (\rho _l^h S_l + \rho _g^h S_g)}{\partial t} + \nabla \cdot (\rho _l^h \vec{q}_l + \rho _g^h \vec{q}_g + j_l^h) = 0, \label{hydrogen_mass_conservation}
\end{align}
where the subscripts $l,g$ denote the liquid and gas phases, and the superscripts $w,h$ denote the water and hydrogen components, respectively. $\phi$ is the porosity, $S_\alpha, \vec{q}_\alpha$ are the saturation and velocity of phase $\alpha$, respectively; $\rho _l^h$ is the dissolved hydrogen mass concentration in the liquid phase; and $j_l^h$ is the diffusion flux of hydrogen in the liquid phase. The Darcy's velocity $\vec{q}_\alpha$ follows the Darcy-Muskat law:
\begin{align}
\vec{q}_\alpha = -K \lambda _\alpha \nabla (P_\alpha - \rho _\alpha \vec{g}), \hspace{5mm} \alpha = l,g ,
\end{align}
where $K$ is the absolute permeability, $\lambda _\alpha$, $P_\alpha$, and $\rho _\alpha$ are the mobility, pressure, and density of phase $\alpha$, and $\vec{g}$ is the gravitational acceleration. The mobility $\lambda _\alpha$ of phase $\alpha$ is defined as the ratio between the phase relative permeability $k_{r\alpha}$ and the phase viscosity $\mu _\alpha$: $\lambda _\alpha = k_{r\alpha}/\mu _\alpha$. Using Fick's law, the diffusion flux of hydrogen in liquid $j^h_l$ in equations \cref{water_mass_conservation,hydrogen_mass_conservation} can be expressed as:
\begin{align}
j_l^h = -\phi S_l D^h_l \nabla \rho ^h_l,
\end{align}
where $D_l^h$ is hydrogen molecular diffusion coefficient in liquid. Since we assume incompressibility of the liquid phase, the mass density of the water component in the liquid phase is constant, i.e. $\rho _l^w = \rho _w^{std}$. To capture capillarity effect, the jump in the pressure at the interface of the two phases is modeled by the relation:
\begin{align}
&P_g = P_l + P_c(S_l).
\end{align} 
where $P_c$ is the capillary pressure. Additionally, we have the constraints 
\begin{align}
S_l + S_g = 1.
\end{align}
To close the model, we also need a set of equations for the thermodynamic equilibrium. Neglecting water vapor and assuming low solubility of hydrogen in the liquid phase, Henry's law can be used to connect the gas pressure $P_g$ and the dissolved hydrogen mass concentration in liquid $\rho _l^h$:
\begin{align}
\rho _l^h = C_h P_g, \label{henrys_law}
\end{align} 
where $C_h = HM^h = \rho _w^{std}M^h/M^w$, $H$ is the Henry's law constant, and $M^i, \; i \in \{w,h\}$ are the molar mass of the $i$-th component. Again, since we ignore water vapor in the gas phase, the ideal gas law reads:
\begin{align}
\rho _g^h = \rho _g = C_v P_g,
\end{align}
where $C_v$ is a constant and $C_v = M^h/(RT)$; T is the temperature and R the ideal gas constant.

\subsection{Nonlinear Complementarity Problem}
To handle phase transitions, we introduce the following nonlinear complementarity \textcolor{black}{constraint problem}
\begin{align}
C_h P_g - \rho _l^h \geq 0, \hspace{3mm} 1 - S_l \geq 0, \hspace{3mm} (1 - S_l)(C_h P_g - \rho _l^h) = 0. \label{ncp}
\end{align}
Equivalently, we can rewrite the above equation using the min function as in \cite{BenGharbia14,Lauser11}
\begin{align}
\min (1 - S_l, C_h P_g - \rho _l^h) = 0. \label{cfun}
\end{align}
Although one can use other types of complementarity functions, the min function is convenient because of its piece-wise linearity with respect to the variable $S_l$ and $\rho _h^l$, which simplifies the computation of the Jacobian in each nonlinear iteration. When the gas phase is not present, we have $C_h P_g - \rho _l^h > 0$ since $\rho_l^h = 0$, and equation \cref{cfun} reduces to $1 - S_l = 0$. When the gas phase appears, $1 - S_l > 0$ and the constraint equation is governed by Henry's law \cref{henrys_law}.

\subsection{Relative Permeability Curves}
In this paper, we use two different models for relative permeability terms.
\begin{itemize}
\item Power law
\begin{align}
&k_{rl} = S_{le}^2, \hspace{5mm} k_{rg} = (1-S_{le})^2, \label{power_law} \\
&S_{le} = \dfrac{S_l - S_{lr}}{1 - S_{lr} - S_{gr}}.
\end{align}
\item Van Genuchten law \cite{VanGenuchten80}
\begin{align}
&k_{rl} = \sqrt{S_{le}}\Big(1 - \big(1 - S_{le}^{1/m}\big)^m\Big)^2, \\
&k_{rg} = \sqrt{1 - S_{le}}\Big( 1- S_{le}^{1/m} \Big)^{2m}, \\
&m = 1 - \dfrac{1}{n},
\end{align}
\end{itemize}
where $S_{le}$ is the effective liquid saturation, and $S_{lr}, S_{gr} \in [0,1]$ are the residual saturations of the liquid and gas phase, respectively.

\subsection{Capillary Pressure}
We employ two models for capillary pressure
\begin{itemize}
\item Linear model
\begin{align}
P_c = P_r(1 - S_{le}).
\end{align}
\item Van Genuchten model \cite{VanGenuchten80}
\begin{align}
P_c = P_r\Big(S_{le}^{-1/m} - 1\Big)^{1/n}.
\end{align}
\end{itemize}
\textcolor{black}{where $P_r$ is the entry pressure.} Notice that the function $P_c(S_l)$ in the Van Genuchten model is only defined for $S_l \in [0,1]$ and $P_c^\prime$ is unbounded near 0 and 1. Thus, it is necessary to modify the model to limit the growth of $P_c^\prime$ and extend it for $S_l \in \mathbb{R}$, since the value of $S_l$ can become larger than 1 or less than 0 during the nonlinear iteration. We used a regularization as presented in \cite{Marchand14} with parameter $\epsilon = 10^{-5}$.
%\begin{itemize}
%\item For $S_l \in [S_{gr}, 1 - S_{lr}]$
%\begin{align}
%& \tilde{S} := S_{gr} + (1 - \epsilon)(S_l - S_{gr}) + \dfrac{\epsilon}{2}(1 - S_{gr} - S_{lr}) \\
%& \tilde{P}_c(S_l) = P_c(\tilde{S}) - P_c\Big( S_gr + \dfrac{\epsilon}{2}(1 - S_{gr} - S_{lr}) \Big)
%\end{align}
%\item For $S_l < S_{gr}$
%\begin{align}
%\tilde{P}_c(S) = \tilde{P}_c(S_{gr}) + \tilde{P}_c^\prime (S_{gr})(S_l - S_{gr})
%\end{align}
%\item For $S_l > 1 - S_{lr}$
%\begin{align}
%\tilde{P}_c(S) = \tilde{P}_c(1 - S_{lr}) + \tilde{P}_c^\prime (1 - S_{lr})(S_l - 1 + S_{lr})
%\end{align}
%\end{itemize}

\subsection{Primary Variables}
There are many ways to choose a set of primary variables, depending on the problem formulation and applications. In our model example, a convenient choice is the liquid saturation, liquid pressure, and the concentration of hydrogen in the liquid phase. We have our solution vector $u = \{P_l, S_l, \rho _l^h\}$. Unlike in other methods such as primary variable switching, in the NCP approach, the choice of primary variable is fixed throughout the simulation. This is an important feature for success of our multilevel algorithm discussed in section \ref{sec:MGR-apply}. 

\section{Solution Algorithm}
\label{sec:algorithm}
In this paper, we consider solving the coupled system consisting of \cref{water_mass_conservation,hydrogen_mass_conservation,ncp} fully implicitly. We use a cell-centered finite volume method for spatial discretization, as it is a natural way to preserve the mass conservation property of the balance equations \cref{water_mass_conservation,hydrogen_mass_conservation}. In addition, it can deal with the case of discontinuous permeability coefficients, and it is relatively straightforward to implement. 
%Under appropriate assumptions, this method also falls into the mixed finite element framework \cite{Peaceman77,Russell83}. 
For the time domain, we employ the backward Euler method to avoid the CFL stability restriction of the time step.
\subsection{Semi-smooth Newton's Method}
We want to solve the system
\begin{displaymath}
\vec{R}(u) = \left \{
\begin{array}{lr}
\vec{H}(u) = 0 & \text{(from the PDEs)} \\
\boldsymbol{\Theta}(u) = \min (\vec{F}, \vec{G}) = 0 & \text{(from the constraints)}
\end{array}
\right.
\end{displaymath}
in which $\vec{F}$ and $\vec{G}$ are discrete functions of $1 - S_l$ and $C_g P_g - \rho _l^h$ respectively, and 
$\vec{R}(u)$ is the residual function. A straightforward approach for solving nonlinear systems of equations is the Newton's method, which requires solution of a linear system at each iteration $k$:
\begin{align}
\dfrac{\partial \vec{R}}{\partial u}\Big|_{u=u_k} \delta u = - \vec{R}(u_k) .
\end{align}
This method requires that the Jacobian $\partial \vec{R}/\partial u$ is defined everywhere. In the case of NCP formulation, the constraints $\boldsymbol{\Theta}$ are only differentiable almost everywhere, and we will need to consider a \textit{semi-smooth} Newton's method instead. The procedure for the semi-smooth Newton's method is similar to that for Newton's method, except that we substitute the derivative $\boldsymbol{\Theta}^\prime$ with the subdifferential $\partial \boldsymbol{\Theta}$ when the function $\boldsymbol{\Theta}$ is non-differentiable. Let $F: \mathbb{R}^n \rightarrow \mathbb{R}^n$ be a locally Lipschitz-continuous function and $D_F$ be the set where $F$ is differentiable; the B-subdifferential of $F$ at $x$ is defined as the set
\begin{align*}
\partial _B F(x) := \{G \in \mathbb{R}^{n\times n}: \exists \; x_k \in D_F \text{ with } x_k \rightarrow x, \nabla F(x_k) \rightarrow G \}\;.
\end{align*}
Below is the algorithm for the semi-smooth Newton's method as described in \cite{BenGharbia14}.\\
\begin{algorithm}[H]
\TitleOfAlgo{Semi-smooth Newton Method}
\While{k $<$ max\_iter and res $> \epsilon$}{
  (1) Define the index sets $A^k$ and $I^k$: \\
  \hspace{6mm}$A^k := \lbrace j: \vec{F}_j(u^k) \ge \vec{G}_j(u^k) \rbrace$,  $I^k := \lbrace j: \vec{F}_j(u^k) < \vec{G}_j(u^k) \rbrace$ \vspace{1mm} \\
  (2) Select an element $J^k \in \partial \boldsymbol\Theta (u^k)$ such that its $j$th \textcolor{black}{row} is equal to\\
  \hspace{6mm}$\vec{F}^\prime _j(u^k)$ if $j \in I^k$, $\vec{G}^\prime _j(u^k)$ if $j \in A^k$ \vspace{1mm}\\
  (3) Solve the system \\
     \hspace{6mm}$\vec{H}^\prime (u^k) \bigtriangleup u^k = - \vec{H}(u^k)$ \\
     \hspace{6mm}$J^k \bigtriangleup u^k = - \boldsymbol\Theta (u^k)$ \vspace{1mm}\\
  (4) Update $u^{k+1}$ \\
     \hspace{6mm}$u^{k+1} = u^k + \bigtriangleup u^k$
}
\end{algorithm}
\noindent For our two-phase, two-component model, the active set $A^k$ corresponds to the set of \textcolor{black}{last rows} where the gas phase is present. The general semi-smooth Newton's method has superlinear convergence for semi-smooth functions, and quadratic for strongly semi-smooth functions. Definitions of semi-smooth and strongly semi-smooth are given in \cite{Qi99}, and a complete treatment of the semi-smooth Newton's method with active set strategy is presented in \cite{Hintermuller02}. 
\par The linear system resulting from taking the subdifferential  $\partial \vec{R}/\partial u$ is often very difficult to solve using iterative methods, and preconditioning is critical for rapid convergence of Krylov subspace methods such as GMRES. In the next section, we discuss the linear system arising from the semi-smooth Newton's method and give a detailed description of the solution algorithms we will use to solve this system.

\subsection{Linear System}
\par Assuming that each physical variable is ordered lexicograhpically, then each nonlinear iteration entails the solution of a discrete version of a block linear system of the form
\begin{align}
\begin{pmatrix}
A_{11} & A_{12} & A_{13} \\
A_{21} & A_{22} & A_{23} \\
A_{31} & A_{32} & A_{33}
\end{pmatrix} 
\begin{pmatrix}
u_1 \\
u_2 \\
u_3
\end{pmatrix} = \begin{pmatrix}
f_1 \\
f_2 \\
f_3
\end{pmatrix}, \label{linear_system}
\end{align}
in which the matrices in the first two rows are the discretized version of the linearized operators from the PDEs. Let $\delta P_l, \delta S_l, \delta \rho _l^h$ be the updates for pressure, saturation, and hydrogen density at each nonlinear step. Using Taylor expansion and keeping only the linear terms, we have
\begin{align*}
&A_{11} = - \nabla \cdot (\rho _l^w \bm{K} \lambda _l \nabla \mathcolor{black}{\delta P_l}), \\
&A_{12} = \phi \dfrac{\partial}{\partial t} (\rho ^w_l \mathcolor{black}{\delta S_l}) - \nabla \cdot (\rho ^w_l \bm{K} \lambda^\prime_l \nabla \tilde{P}_l \mathcolor{black}{\delta S_l}) + \nabla \cdot (\phi D^h_l \nabla \tilde{\rho} ^h_l \mathcolor{black}{\delta S_l} ), \\
&A_{13} = \nabla \cdot (\phi S_l D^h_l \nabla \mathcolor{black}{\delta \rho ^h_l}), \\
&A_{21} = \phi \dfrac{\partial}{\partial t} (S_g C_g \mathcolor{black}{\delta P_l}) - \nabla \cdot (\rho ^h_l \bm{K} \lambda _l \nabla \mathcolor{black}{\delta P_l}) - \nabla \cdot (\rho ^h_g \bm{K} \lambda _g \nabla \mathcolor{black}{\delta P_l}) - \nabla \cdot (C_g \bm{K} \lambda _g \nabla P_g \mathcolor{black}{\delta P_l}), \\
&A_{22} = \phi \dfrac{\partial}{\partial t}((\rho ^h_l - \rho ^h_g)\mathcolor{black}{\delta S_l}) - \nabla \cdot ( \rho ^h_l \bm{K} \lambda ^\prime _l \nabla P_l \mathcolor{black}{\delta S_l}) - \nabla \cdot (\rho ^h_g \bm{K} \lambda ^\prime _g \nabla P_g \mathcolor{black}{\delta S_l}) \\
&- \nabla \cdot (C_g P_c^\prime \bm{K} \lambda _g \nabla P_g \mathcolor{black}{\delta S_l}) - \nabla \cdot (\rho ^h_g \bm{K} \lambda _g \nabla P_c^\prime \mathcolor{black}{\delta S_l}) - \nabla \cdot (\rho ^h_g \bm{K} \lambda _g P_c^\prime \nabla \mathcolor{black}{\delta S_l}) \\
&- \nabla \cdot (\phi D^h_l \nabla \rho ^h_l \mathcolor{black}{\delta S_l}), \\
&A_{23} = \phi \dfrac{\partial}{\partial t} (\tilde{S}_l \mathcolor{black}{\delta x^h_l}) - \nabla \cdot (\phi \tilde{S}_l D^h_l \nabla \mathcolor{black}{\delta \rho ^h_l}).
\end{align*}
All the coefficients in the above equations are evaluated at the linearization point $\{\tilde{P}_l, \tilde{S}_l, \tilde{\rho}_l^h\}$. From these operators, we can make some important observations:
\begin{itemize}
\item The global matrix is non-symmetric and indefinite.
\item The block $A_{11}$ has the structure of a discrete purely elliptic problem for pressure.
\item The coupling block $A_{12}$ has the structure of a discrete first-order hyperbolic problem in the liquid phase saturation.
\item The coupling block $A_{21}$ has the structure of a discrete parabolic problem in the wetting phase pressure.
\item The block $A_{22}$ has the structure of a discrete parabolic (convection-diffusion) problem for saturation when capillary pressure is a non-constant function of the saturation. When capillary pressure is zero or a constant, $P_c^\prime = 0$ and there is no diffusion term, the block has the form of a hyperbolic problem.
\item \textcolor{black}{The entries of the blocks with respect to the dissolved hydrogen mass density $A_{13}, A_{23}$ are small with respect to the diagonal block $A_{11}$}, and only play a significant role in the regions where the gas phase does not exist.
\end{itemize}
These observations will help us motivate the development of our new method in the next section.
\par Besides the blocks associated with the PDEs, we also need to consider those in last row of the matrix in \cref{linear_system}, which are derived from the discrete version of the complementarity constraint equation \cref{ncp}. When the gas phase does not exist, we have
\begin{align*}
&A_{31} = 0, \hspace{5mm} A_{32} = - \delta S_l, \hspace{5mm} A_{33} = 0,
\end{align*}
and when the gas phase is present, these blocks become
\begin{align*}
&A_{31} = H \delta S_l, \\
&A_{32} = H P_c^\prime \delta S_l, \\
&A_{33} = \dfrac{M^h \rho _l^w}{M^w} \delta x^h_l.
\end{align*}
In matrix form, the blocks $A_{31}, A_{32}$, and $A_{33}$ are diagonal matrices, since the constraints are local. Again, because a phase can disappear, the block $A_{33}$ is not guaranteed to be non-singular. In fact, when this happens, the rows corresponding to the cells where a phase disappears have zero diagonal values. Thus, we can split the last row into two separate sets: the set with zeros on the diagonal of $A_{33}$ and its complement. Rewriting the matrix $A$ using this splitting we have:
\begin{align}
A = \begin{blockarray}{cccccc}
\begin{block}{(cccc)cc}
A_{11} & A_{12} & A_{13} & A_{14} \\
A_{21} & A_{22} & A_{23} & A_{24} \\
C_{31} & C_{32} & C_{33} & 0\\
C_{41} & C_{42} & 0 & 0\\
\end{block}
\end{blockarray}. \label{2p2c_matrix}
\end{align}
\textcolor{black}{Let $N$ be the number of elements in the mesh and $M$ be the number of cells in which the gas phase is present. Then, the size of the matrix $A$ is $3N \times 3N$. The blocks $A_{ij}$, $i = 1,2$, $j = 1,2$ have the size of $N \times N$. The pressure - hydrogen mass concentration $A_{13}, A_{14}$ and saturation - hydrogen mass concentration $A_{23}, A_{24}$ coupling blocks have the size of $N \times M$ and $N \times (N-M)$, respectively. The hydrogen mass concentration - pressure $C_{31}, C_{41}$ and hydrogen mass concentration - saturation $C_{32}, C_{42}$ constraint blocks have the size of $M \times N$ and $(N-M) \times N$, respectively. The block $C_{33}$ is a diagonal matrix of size $M \times M$ which contains only non-zero diagonal values of $A_{33}$. The $0$ block on the diagonal has the size of $(N-M) \times (N-M)$.} \textcolor{black}{Since $A$ has zeros on its diagonal}, it is clear that we cannot use classical AMG algorithms to solve this system. In the past, since much of the focus was to find a formulation that can take into account all the complex physics involved in simulating compositional multiphase flow with miscibility and phase transitions, there has not been a lot of work in designing optimal preconditioners for this type of linear system. Recently, there have been some development of algebraic multigrid preconditioners such as two-stage preconditioning \cite{Stueben07,Wang17} and block factorization \cite{Bui17} for immiscible two-phase flow. Yet, these methods have not been applied successfully to the problems considered here. The most popular and robust method is still using the incomplete factorization (ILU) of the global matrix $A$ as a preconditioner to GMRES. In this work, we seek to develop a new algebraic multigrid preconditioner based on multigrid reduction for the linear system arising in the semi-smooth Newton's method.

\section{Multigrid Reduction}
\label{sec:MGR}
The idea of multigrid reduction (MGR) has been around for a long time, tracing back to the work of Ries and Trottenberg \cite{Ries79,Ries83}. Recently, it has gained more attention through the work on multigrid reduction in time by Falgout et. al. \cite{Falgout14}. In this section, we summarize the approach for the case of two-level reduction. Given a matrix $A$, we have the C-F splitting 
\begin{align}
A = \begin{pmatrix}
A_{ff} & A_{fc} \\
A_{cf} & A_{cc}
\end{pmatrix} = \begin{pmatrix}
I_{ff} &  0\\
A_{cf}A_{ff}^{-1} & I_{cc}
\end{pmatrix} \begin{pmatrix}
A_{ff} & 0 \\
0 & S
\end{pmatrix} \begin{pmatrix}
I_{ff} & A_{ff}^{-1}A_{fc}  \\
0& I_{cc}
\end{pmatrix},
\end{align}
where $I_{cc}$ and $I_{ff}$ are identity \textcolor{black}{matrices} and $S = A_{cc} - A_{cf} A_{ff}^{-1} A_{fc}$ is the Schur complement.
\par \noindent We can define the ideal interpolation and restriction operators by 
\begin{align}
P = \begin{pmatrix}
-A_{ff}^{-1} A_{fc}\\
I_{cc}
\end{pmatrix}, \hspace{5mm} R = \begin{pmatrix}
-A_{cf}A_{ff}^{-1} &I_{cc}
\end{pmatrix}.
\end{align}
Additionally, define the injection operator as $Q = \begin{pmatrix}
 I_{ff} & 0
\end{pmatrix}^T$. Then since $A_{ff} = Q^TAQ$ and $S = RAP$, it is simple to derive that
\begin{equation*}
A^{-1}=P(RAP)^{-1}R+Q(Q^TAQ)^{-1}Q^T,
\end{equation*}
and 
\begin{align}
  0 = I - A^{-1}A
  &= I - P(RAP)^{-1}RA-Q(Q^TAQ)^{-1}Q^TA\label{eq:MGR-add}\\
  &=(I - P(RAP)^{-1}RA)(I-Q(Q^TAQ)^{-1}Q^TA)\label{eq:MGR-mul1}\\
  &=(I-Q(Q^TAQ)^{-1}Q^TA)(I - P(RAP)^{-1}RA),\label{eq:MGR-mul2}
\end{align}
where the equivalence occurs since $RAQ=Q^TAP=0$. This identity defines the
two-level multigrid method with the ideal Petrov-Galerkin coarse-grid operator
$RAP$ and the F-relaxation $Q(Q^TAQ)^{-1}Q^T$.  Equation \cref{eq:MGR-add}
is the additive MGR identity and \cref{eq:MGR-mul1} and \cref{eq:MGR-mul2} are
multiplicative identities with pre-smoothing and post-smoothing. However, constructing ideal interpolation and restriction operators is impractical, and we need to approximate these operators.  In practice, MGR methods replace ideal restriction and prolongation with approximations $R$ and $P$ respectively, where
\begin{equation}
\label{eq:RPoperator1}
P= 
\begin{pmatrix}
W_{p}
\\ I_{cc},
\end{pmatrix},\quad
R =
\begin{pmatrix}
W_{r} & I_{cc}
\end{pmatrix}.
\end{equation}
There are many ways to construct the restriction $R$ and interpolation $P$ operators. Here, we have only experimented with two options
\begin{align}
W_r = 0, \hskip2ex W_p = -D_{ff}^{-1} A_{fc}, \label{injective_rp}
\end{align}
and 
\begin{align}
W_r = -A_{cf}D_{ff}^{-1}, \hskip2ex W_p = -D_{ff}^{-1} A_{fc}, \label{non_injective_rp}
\end{align}
where $D_{ff} = diag(A_{ff})$.

The F-relaxation in \cref{eq:MGR-mul1} and \cref{eq:MGR-mul2} is also generally replaced with a more efficient method and often extended to all unknowns, not just F-points. We can solve with (block) Jacobi, (block) Gauss-Seidel, ILU, or AMG.

{\color{black} The coarse grid operator $A_{c} = RAP$ could also be considered as an approximation to the Schur complement. 
There are many proposed approximations to Schur complement, including several based on multigrid ideas \cite{bank1999incomplete,chow2003multilevel,reusken1996multigrid,reusken1999approximate,rusten1992preconditioned,Wagner1997}. 
The physics-based approximation for the Schur complement are of interest general settings. There are lots literatures for saddle point problems  \cite{benzi2008some,elman2006block,elman2008taxonomy,sloan1986algorithm}. Another interesting direction is based on the work on Block Factorized Sparse Approximate Inverse (Block FSAI) preconditioners in \cite{Ferronato14}.  In the context of MGR, the Block FSAI could be used to construct good approximation to the restriction and interpolation operators.}

In general, we define the MGR operator in either pre-smoothing or post-smoothing
form by
%\begin{equation}
%  \label{eq:MGR1}
%  (M_{MGR}^1)^{-1}A = PM_{C}^{-1}RA+SM_{n}^{-1}S^TA
%\end{equation}
%or
\begin{align}
  I - M_{MGR}^{-1}A &= (I - PM_{c}^{-1}RA)(I-M_{f}^{-1}A) , \label{eq:MGR-pre}\\
  I - M_{MGR}^{-1}A &= (I-M_{f}^{-1}A)(I - PM_{c}^{-1}RA) , \label{eq:MGR-post}
\end{align}
where $M_{c}^{-1}\approx (RAP)^{-1}$ is the coarse grid correction and
$M_{f}^{-1}$ is the smoother. The two grid solve consists of an F-relaxation followed by a coarse-grid correction can be presented as follows

\begin{algorithm}[H]
\TitleOfAlgo{MGR preconditioner with presmoothing}
Let $r=b$ and $e_{0} = 0$\;
  (1) {\bf Global Relaxation:} \\
  \hspace{6mm}$e_{1}\leftarrow e_{0} + M^{-1}r$, where $M=\text{blockdiag}(A)$ \\
  (2) {\bf F-Relaxation:} \\
  \hspace{6mm}$e_{2}\leftarrow e_{1} + QM_{ff}^{-1}Q^{T}(r-Ae_{1})$ where $M_{ff} = \text{blockdiag}(A_{ff})$ \vspace{1mm}\\
  (3) {\bf Coarse Grid Correction:}\\
     \hspace{6mm}$e_{3}\leftarrow e_{2} + PM_{c}^{-1}R_{I}(r-Ae_{2})$ where $M_{c}^{-1}$ is approximated by a classical AMG method \\
\end{algorithm}

%\begin{algorithm}[H]
%  Let $r=b$ and $e_{0} = 0$\;\\
%  \;\\
%  {\bf F-Relaxation:} $e_{2}\leftarrow e_{1} + SM_{nn}^{-1}S^{T}(r-Ae_{1})$ where $M_{nn} = \text{blockdiag}(A_{nn})$\;\\
%  {\bf Coarse Grid Correction:} $e_{3}\leftarrow e_{2} + PM_{C}^{-1}R_{I}(r-Ae_{2})$ where $M_{C}^{-1}$ is approximated by a classical AMG method\;\\
%\end{algorithm}

\textcolor{black}{The appeal of the MGR approach is that it provides a flexible framework for choosing the coarse/fine grids, the interpolation and restriction operators, and the solver for the coarse/fine grids. For example, if one chooses to extend the F-points to all unknowns (rather than the complement of the C-points), $W_p = 0$, $W_r = 0$ for the interpolation and restriction operators, and ILU0 and AMG for the F-relaxation and coarse-grid solve, respectively, then the MGR method is equivalent to the CPR-AMG approach. The block factorization method in \cite{Bui17} is another variant of the MGR approach, which uses the C-points for the pressure and F-points for the saturation unknowns and $W_r = - A_{cf}D_{ff}^{-1}$, $W_p = 0$. Another advantage of the MGR approach is that it is an algebraic method and unlike geometric multigrid, it can be used as a ``black-box'' solver for general geometries and grid types.}

\section{MGR for the two-phase, two-component model}
\label{sec:MGR-apply}
\begin{figure}[H]
\centering
\includegraphics[width=0.5\textwidth]{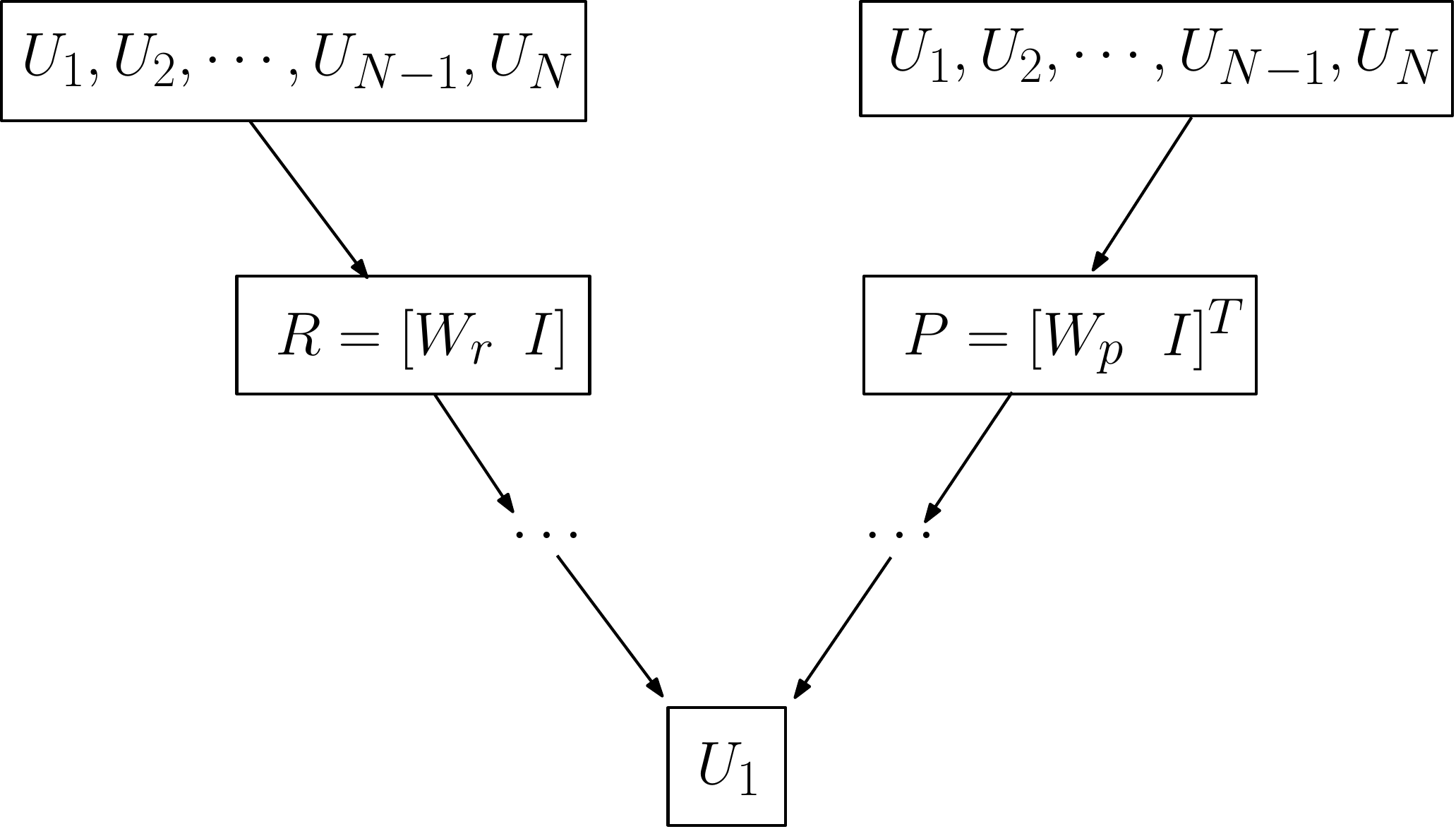}
\captionsetup{justification=centering}
\caption{Multigrid reduction V-cycle} \label{MGR_vcycle}
\end{figure}

Again, we need to solve the linear system $A\mathbf{u} = \mathbf{f}$, in which the matrix $A$ is given in \cref{2p2c_matrix}.
The first step of multigrid reduction aims to eliminate the third row, corresponding to the constraints in the cells where all the phases exist. Thus, we have the following splitting:
\begin{align}
A = \begin{blockarray}{cccccc}
\begin{block}{(cccc)cc}
A_{11} & A_{12} & A_{13} & A_{14} & & C\\
A_{21} & A_{22} & A_{23} & A_{24} & & C\\
C_{31} & C_{32} & C_{33} & 0 & & F\\
C_{41} & C_{42} & 0 & 0 & & C\\
\end{block}
\end{blockarray}
\end{align}
Note that the last column indicates the C-F splitting we use for this case. The Schur complement after the reduction step reads
\begin{align}
S^1 &= RAP =  \begin{pmatrix}
A_{11} & A_{12} & A_{14} \\
A_{21} & A_{22} & A_{24} \\
C_{41} & C_{42} & 0 \\
\end{pmatrix} - \begin{pmatrix}
A_{13} \\
A_{23} \\
0
\end{pmatrix} C_{33}^{-1} \begin{pmatrix}
C_{31} & C_{32} & 0
\end{pmatrix} \\
& = \begin{pmatrix}
A_{11} - A_{13}C_{33}^{-1}C_{31} & A_{12} - A_{13}C_{33}^{-1}C_{32} & A_{14} \\
A_{21} - A_{23}C_{33}^{-1}C_{31} & A_{22} - A_{23}C_{33}^{-1}C_{32} & A_{24} \\
C_{41} & C_{42} & 0
\end{pmatrix}.
\end{align}
%\begin{align}
%&S_{11} = A_{11} - A_{13}C_{33}^{-1}C_{31} \\
%&S_{12} = A_{12} - A_{13}C_{33}^{-1}C_{32} \\
%&S_{13} = A_{14} \\
%&S_{21} = A_{21} - A_{23}C_{33}^{-1}C_{31} \\
%&S_{22} = A_{22} - A_{23}C_{33}^{-1}C_{32} \\
%&S_{23} = A_{24} \\
%&S_{31} = C_{41} \\
%&S_{32} = C_{42} \\
%&S_{33} = 0
%\end{align}
Again, the operators $R$ and $P$ come from equations \cref{injective_rp,non_injective_rp}. Note that this reduction step is exact since $C_{33}$ is a diagonal matrix. However, we have not eliminated the zero diagonal values after the first reduction step. Next, we eliminate the saturation block with the following C-F splitting:
\begin{align}
S^1 = \begin{blockarray}{ccccc}
\begin{block}{(ccc)cc}
S_{11} & S_{12} & A_{14} & & C\\
S_{21} & S_{22} & A_{24} & & F\\
C_{41} & C_{42} & 0 & & C\\
\end{block}
\end{blockarray}
\end{align}
The Schur complement \textcolor{black}{(also the coarse grid)} for the second level of multigrid reduction reads
\begin{align}
S^2 &= RS^1P = \begin{pmatrix}
S_{11} & A_{14} \\
C_{41} & 0
\end{pmatrix} - \begin{pmatrix}
S_{12} \\
C_{42}
\end{pmatrix} {\tilde{S}_{22}^{-1}} \begin{pmatrix}
S_{21} & A_{24}
\end{pmatrix} \\
& = \begin{pmatrix}
S_{11} - S_{12} \tilde{S}_{22}^{-1}S_{21} & A_{14} - S_{12} \tilde{S}_{22}^{-1}A_{24} \\
C_{41} - C_{42} \tilde{S}_{22}^{-1}S_{21} & -C_{42} \tilde{S}_{22}^{-1}A_{24}
\end{pmatrix},
\end{align}
\textcolor{black}{where $\tilde{S}_{22}^{-1}$ \textcolor{black}{is some approximation} of $S_{22}^{-1}$ \textcolor{black}{to compute} $R$ and $P$ \textcolor{black}{from (\ref{eq:RPoperator1})}. In the F-relaxation step, \textcolor{black}{the action of the} saturation block $S_{22}^{-1}$ is \textcolor{black}{achieved by} one V-cycle of AMG.} In the equation above, the presence of the constraints in the matrix $S^2$ makes it non-elliptic, and therefore, we cannot solve it using AMG, although it no longer has zeros on the diagonal. The final reduction step is employed to eliminate these constraints by putting them as F-points:
\begin{align}
S^2 &= 
\begin{blockarray}{cccc}
\begin{block}{(cc)cc}
S^2_{11} & S^2_{12} & & C\\
S^2_{21} & S^2_{22} & & F\\
\end{block}
\end{blockarray} \\
S^3 &= RS^2P = S^2_{11} - S^2_{12}(\tilde{S}^2_{22})^{-1}S^2_{21}
\end{align}
The Schur complement at the last level $S^3$ can be solved using AMG. A \textcolor{black}{schema} of a multi-level reduction approach is illustrated in Fig. \ref{MGR_vcycle}.
\noindent At each level, we reduce one or more variables until we obtain the variable associated with the elliptic operator, which is the pressure variable in our case. This means that in Fig. \ref{MGR_vcycle}, $U_1 \equiv P_l$.

\section{Numerical Experiments}
\label{sec:result}
In this section, we perform numerical experiments to show the efficiency of the multigrid reduction approach (hereinafter referred as hypreMGR) in solving the linear equations arising in various flow scenarios. The algorithm is implemented as a part of Hypre \cite{Falgout06,Falgout02}. For all AMG solve steps, we use BoomerAMG \cite{Henson00}, also included in Hypre. The two-phase, two-component flow using the NCP approach is implemented in Amanzi, a parallel open-source multi-physics C++ code developed as a part of the ASCEM project \cite{Amanzi}. Although Amanzi was first designed for simulation of subsurface flow and reactive transport, its modular framework and concept of process kernels \cite{ECoon_JDMoulton_SPainter_2014a} allow new physics to be added relatively easily for other applications. The simulator employed in this work is one such example. Amanzi works on a variety of platforms, from laptops to supercomputers. It also leverages several popular packages for mesh infrastructure and solvers through a unified input file. \textcolor{black}{Again, due to the presence of zeros on the diagonal of the generalized Jacobian matrix $A$, we cannot use current AMG solvers such as Hypre's BoomerAMG and Trillinos ML for our problem. Even when we eliminate the zeros on the diagonal of $A$ to obtain a smaller system, using HYPRE's BoomerAMG and Trilinos ML for system as a preconditioner, GMRES still fails to converge within 400 iterations. We also exclude ILUt variants (also implemented in Euclid and Trillinos) as they are not robust for the problems considered in this work. Thus, in all of our experiments we compare hypreMGR with the ILU($k$) method from Euclid, which is also a part of Hypre.} ILU($k$) is used sequentially for all the examples. We experiment with different levels of fill $k$ and report the results for the minimum $k$ that is sufficient for GMRES to converge within 400 iterations throughout the simulation in each test case. \textcolor{black}{We also try ILU2 (see \cite{Konshin15}), which is a two-parameter modification of ILUt designed for nonsymmetric saddle-point problems, and find its performance comparable to that of ILU($k$) presented in this section with appropriate thresholds. Yet, it is not clear how to choose the optimal parameters of ILU2 for the problems presented here, since they are dependent on the characteristics of the problem, i.e. advection-dominated or diffusion-dominated, and also on the mesh size. Due to these complexities, we believe a fair comparison between ILU2 and our method is beyond the scope of this paper.} GMRES is provided within Amanzi. For simplicity, we employ structured Cartesian grids for the test cases, but we can also use unstructured K-orthogonal grids. For parallel results, the test cases are run on Syrah, a Cray system with 5,184 Intel Xeon E5-2670 cores at the Lawrence Livermore National Laboratory Computing Center. Amanzi and other libraries are compiled with OpenMPI 1.6.5 and gcc-4.9.2. The total time is measured in seconds. 
\par This section has four parts. In the first part, we show the results for an unsaturated flow problem with no phase appearance/disappearance. In the second part, we report the results for the saturated flow problem in which the gas phase appears by injection and then disappears after the injection is stopped. These two test cases were originally presented in the MoMaS gas benchmark project \cite{BourgeatBenchmark13}. The third example includes a three-dimensional problem with parameters generally used in reservoir simulation. In the last part, we examine the scalability of the multigrid reduction approach.
\par Unless specified otherwise, for all of the simulations presented here, the convergence tolerance for semi-smooth Newton's method is $||F(x)|| \le 10^{-5}$, and the linear tolerance for GMRES is $||J\delta u_k - F(u_k)|| \le 10^{-12}||F(u_k)||$, which is the default in Amanzi. For AMG solves, we use the default parameters in BoomerAMG. The coarsening strategy is the parallel Cleary-Luby-Jones-Plassman (CLJP) coarsening \cite{Cleary98}. The interpolation method is the classical interpolation defined in \cite{RugeStueben86}, and the smoother is the forward hybrid Gauss-Seidel/SOR scheme. The number of V-cycle steps is set to 1.

\subsection{Unsaturated flow}
This test shows a \textcolor{black}{two-dimensional} case in which the water and gas system is initially out of equilibrium and then evolves towards equilibrium. There is no flow in and out of the domain, and there is no phase appearance/disappearance. The detailed set up of the experiment is described below.
\begin{figure}[H]
\centering
\includegraphics[width=0.5\textwidth]{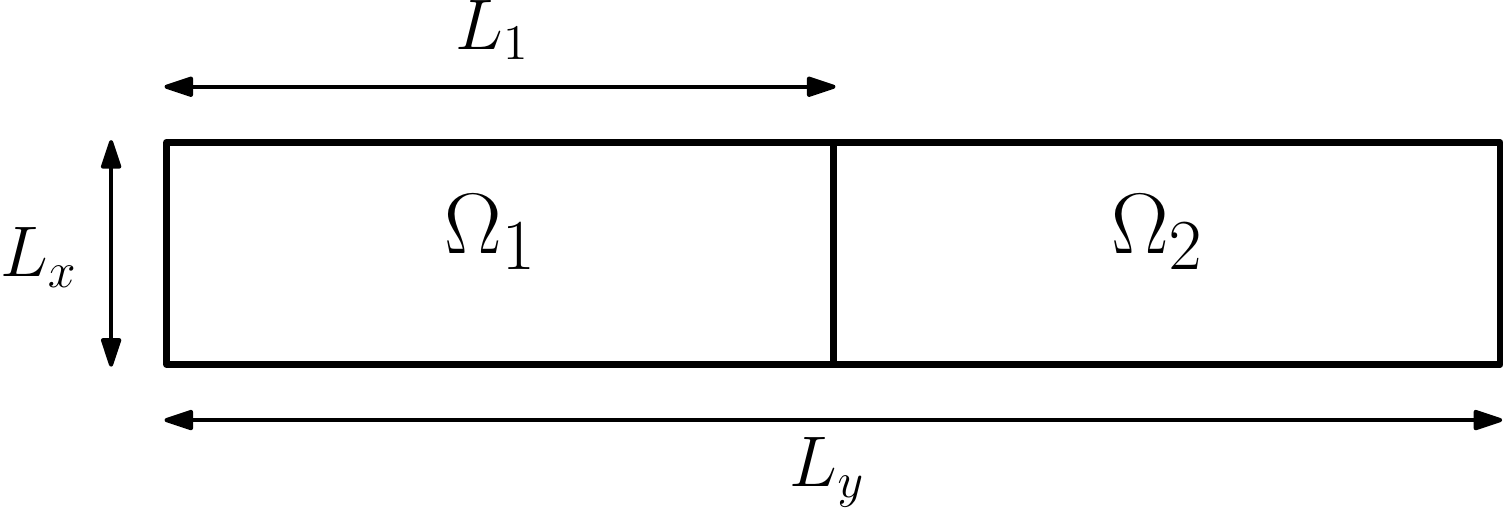}
\captionsetup{justification=centering}
\caption{Porous domain $\Omega$ with two sub-domains $\Omega _1$ and $\Omega _2$.}
\end{figure}
For boundary conditions, we impose no flow condition on the boundary of the whole domain. Denoting $\psi ^w = \rho ^w_l K \lambda _l \nabla P_l + j^h_l$ and $\psi ^h = \rho ^h_l K \lambda _l \nabla P_l + \rho ^h_h K \lambda _g \nabla P_g - j^h_l$, we have $
\psi ^k \cdot \nu = 0, \hspace{2mm} k \in \{w,h\}$ on $\Gamma$. Initial conditions are uniformly constant on each sub-domain $\Omega _1$ and $\Omega _2$:
\begin{itemize}
\item $P_l = P_{l,1}$ and $P_g = P_{g,1}$ on $\Omega _1$.
\item $P_l = P_{l,2} = P_{l,1}$ and $P_g = P_{g,2} \neq P_{g,1}$ on $\Omega _2$.
\end{itemize}
For capillary pressure, we use the Van Genuchten model with $P_r = 2\times 10^6 \text{   Pa}$, $n = 1.54$, $S_{lr} = 0.01$, and $S_{gr} = 0$. The rest of the parameter values are shown in \cref{ini_cond_radon2,param_list_radon2}. We run 5 time steps of size $dT = 10$ seconds. The results are summarized in Table \ref{performance_radon2}. NS denotes the number of nonlinear iterations, LS the number of linear iterations, and LS/NS the average number of linear iterations per nonlinear iterations.
\begin{table}[H]
\begin{minipage}[t]{0.45\textwidth}
\caption{Initial conditions} \label{ini_cond_radon2}
\centering
{\renewcommand{\arraystretch}{1.2}
\begin{tabularx}{0.7\textwidth}{ l X }
  \hline
  $L_1$ & 0.5 m \\
  $L_x$ & 1 m \\
  $L_y$ & 0.1 m \\
  $P_{l,1}$ & $10^{6}$ Pa\\
  $P_{g,1}$ & $1.5 \times 10^6$ Pa\\
  $S_{l,1}$ & $96.2$ \% \\
  $P_{l,2}$ & $10^{6}$ Pa\\
  $P_{g,2}$ & $2.5 \times 10^{6}$ Pa\\
  $S_{l,2}$ & $84.2$ \% \\
  \hline
\end{tabularx}}
\end{minipage}
\begin{minipage}[t]{0.45\textwidth}
\caption{Parameter Values} \label{param_list_radon2}
\centering
{\renewcommand{\arraystretch}{1.2}
\begin{tabularx}{0.7\textwidth}{ l X }
  \hline
  $\bm{K}$ & $1 \times 10^{-16}$ $m^2$ \\
  $\phi$ & 0.3 \\
  $D^h_l$ & $3 \times 10^{-9}$ $m^2/s$ \\
  $\mu _l$ & $1 \times 10^{-9}$ Pa s\\
  $\mu _g$ & $9 \times 10^{-6}$ Pa s\\
  $H$ & $7.65 \times 10^{-6}$ mol/Pa/$m^3$ \\
  $M^h$ & $2 \times 10^{-3}$ kg/mol\\
  $M^w$ & $1 \times 10^{-2}$ kg/mol\\
  $\rho _l^w$ & $10^3$ kg/$m^3$ \\
  \hline
\end{tabularx}}
\end{minipage}
\end{table}
\begin{table}[H]
\centering
\captionsetup{justification=centering}
\caption{Performance of hypreMGR for different mesh sizes} \label{performance_radon2}
{\renewcommand{\arraystretch}{1.2}
\begin{tabular}{| c | c | c | c | c | c | c |}
\hline
\multirow{2}{*}{Mesh size} & \multicolumn{3}{c}{ILU(0)} & \multicolumn{3}{|c|}{hypreMGR} \\
\cline{2-7} 
& Time (s) & LS & LS/NS & Time (s) & LS & LS/NS \\
\hline
%$200 \times 10$ & 11.3 & 497 & 45.2 & 14.6 & 334 & 30.4 \\
%$400 \times 20$ & 74.5 & 1100 & 84.6 & 71.3 & 417 & 32.1 \\
%$800 \times 40$ & 544.4 & 2259 & 173.8 & 284.1 & 409 & 31.5 \\
%$1600 \times 80$ & - & - & - & 1407.8 & 491 & 37.8 \\
$200 \times 10$ & 11.5 & 555 & 50.5 & 10.8 & 445 & 40.5 \\
$400 \times 20$ & 97.3 & 1283 & 98.7 & 42.2 & 458 & 35.2 \\
$800 \times 40$ & 757.4 & 2479 & 190.7 & 180.3 & 557 & 42.8 \\
$1600 \times 80$ & 5666 & 4321 & 332.4 & 801.8 & 569 & 43.8 \\
\hline
\end{tabular}}
\end{table}
In this experiment, because there is no phase disappearance/appearance, the diagonal of the Jacobian does not have any zero, and ILU(0) can be used as a preconditioner. With respect to hypreMGR, we apply two levels of reduction for this problem, one for the constraints with nonzero diagonal values and one for the saturation block. The approximations for the restriction and interpolation operators from \cref{injective_rp} are used in this case. As the liquid saturation does not go to zero, the effect of capillary pressure is small and that makes the system more advection-dominated. Thus, we do not need use AMG to solve for the saturation correction in the F-relaxation step. Here, we found that using three Gauss-Seidel smoothing steps is sufficient except for the finest grid where AMG with a two-level V(3,3)-cycle is used. For the coarse grid, a single AMG V(2,2)-cycle is applied for all the mesh sizes. Table \ref{performance_radon2} indicates that our new algorithm is more efficient both in terms of run time and number of iterations. It also exhibits near optimal scaling with respect to mesh size. hypreMGR is faster across all the meshes both in terms of the run time and average number of linear iterations. For the mesh size of $800 \times 40$, hypreMGR is twice faster in terms of run time, and takes less than a fifth of the number of iterations of ILU(0). For the largest mesh of $1600 \times 80$, ILU(0) is very inefficient, and hypreMGR is about 7 times faster than ILU(0) in terms of run time and number of iterations.

\subsection{Saturated flow with phase appearance}
This test is devoted to describing gas phase appearance produced by injecting pure hydrogen in a two-dimensional homogeneous porous domain $\Omega$, which was initially 100\% saturated by pure water. The porous domain is a rectangle of size $200\text{m} \times 20\text{m}$ with three types of boundaries : $\Gamma _{\text{in}}$ on the left side is the inflow boundary; $\Gamma _{\text{out}}$ on the right side is the outflow boundary; and $\Gamma _{\text{imp}}$ at the top and bottom is the impervious boundary.
\begin{figure}[H]
\centering
\includegraphics[width=0.7\textwidth]{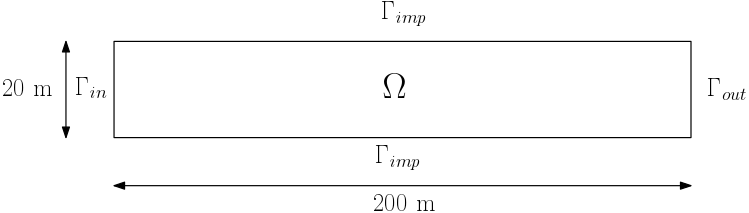}
\captionsetup{justification=centering}
\caption{Core domain for the gas infiltration example.}
\end{figure}
\noindent There is no source terms inside the boundary, and the boundary conditions are as follows
\begin{itemize}
\item No flux on $\Gamma _\text{imp}$
\begin{align}
\psi ^w \cdot \nu = 0 \text{   and   } \psi ^h \cdot \nu = 0
\end{align}
\item Injection of hydrogen on the inlet $\Gamma _{\text{in}}$
\begin{align}
\psi ^w \cdot \nu = 0 \text{   and   } \psi ^h \cdot \nu = 5.57 \times 10^{-6} \hspace{2mm} kg/m^2/year
\end{align}
\item Fixed liquid saturation and pressure on the outlet
\begin{align}
P_l = 10^6 \hspace{2mm} Pa, \hspace{3mm} S_l = 1, \hspace{3mm} \rho ^h_l = 0
\end{align}
\end{itemize}
Initial conditions are uniform throughout the domain, corresponding to a stationary state of saturated liquid and no hydrogen injection
\begin{align}
P_l = 10^6 \hspace{2mm} Pa, \hspace{3mm} S_l = 1, \hspace{3mm} \rho ^h_l = 0
\end{align}
The rest of the physical parameters are given in \cite{BourgeatBenchmark13}.
%The rest of the physical parameters are given in Table \cref{param_list_radon1}
%\begin{table}[H]
%%\begin{minipage}[t]{0.45\textwidth}
%\caption{Parameter Values} \label{param_list_radon1}
%\centering
%{\renewcommand{\arraystretch}{1.2}
%\begin{tabularx}{0.3\textwidth}{ l X }
%  \hline
%  $\bm{K}$ & $5 \times 10^{-20}$ $m^2$ \\
%  $\phi$ & 0.15 \\
%  $D^h_l$ & $3 \times 10^{-9}$ $m^2/s$ \\
%  $\mu _l$ & $1 \times 10^{-9}$ Pa s\\
%  $\mu _g$ & $9 \times 10^{-6}$ Pa s\\
%  $H$ & $7.65 \times 10^{-6}$ mol/Pa/$m^3$ \\
%  $M^h$ & $2 \times 10^{-3}$ kg/mol\\
%  $M^w$ & $1 \times 10^{-2}$ kg/mol\\
%  $\rho _l^w$ & $10^3$ kg/$m^3$ \\
%  \hline
%\end{tabularx}}
%%\end{minipage}
%\end{table}
\par
\begin{figure}[H]
\centering
\includegraphics[width=0.6\textwidth]{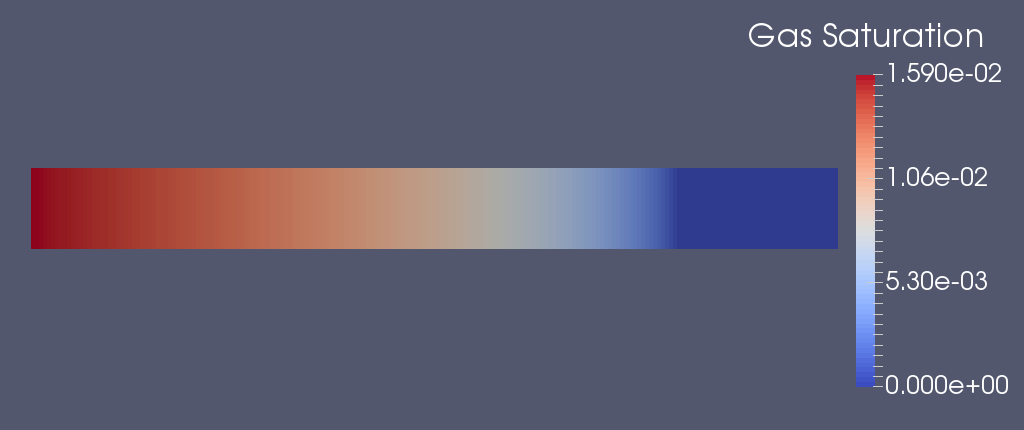}
\captionsetup{justification=centering}
\caption{Gas infiltration after 100 time steps} \label{gas_infiltration}
\centering
\includegraphics[width=0.45\textwidth]{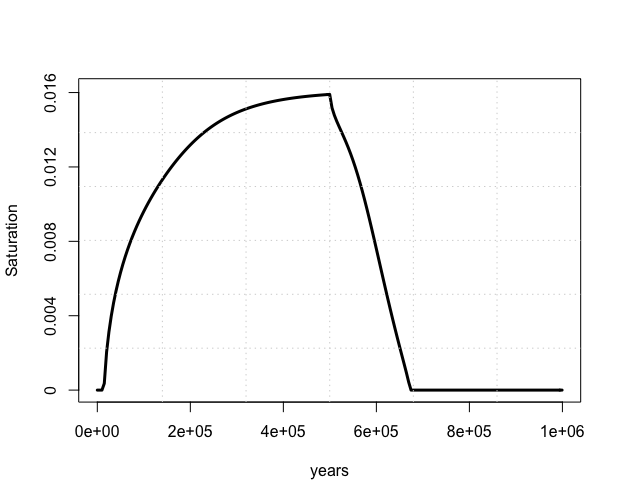}
\includegraphics[width=0.45\textwidth]{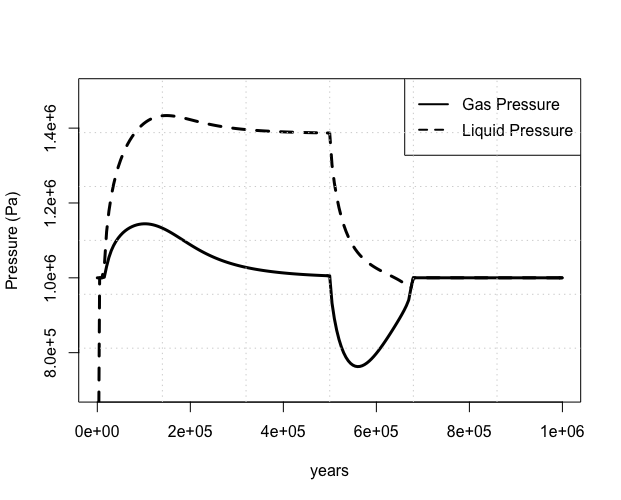}
\caption{Gas saturation and pressure profiles at the first cell over time for the saturated flow case.} \label{radon1_sat_pres_profiles}
\end{figure}

 For the capillary pressure model, we experimented two scenarios: (1) power laws for relative permeabilities as in \cref{power_law} in conjunction with the linear capillary pressure model, and (2) Van Genuchten for both relative permeabilities and capillary pressure model. In the first case, the entry pressure is $P_r = 2 \times 10^6$ Pa, and in the second case, we use $P_r = 2 \times 10^6$ Pa, $n = 1.49$. In both cases, the residual saturations are $S_{gr} = 0$ and $S_{lr} = 0.4$. We run the simulation for 100 time steps of fixed size $dT = 5000$ years. Fig. \ref{gas_infiltration} shows the infiltration of hydrogen after $5 \times 10^5$ years for the second scenario, and the performance of the preconditioners is reported in Table \ref{performance_radon1}. 
\noindent As we can see in Fig. \ref{gas_infiltration}, the left side of the core is infiltrated with hydrogen, while the right side is still fully saturated with water. We also plot the gas saturation and the pressures in the first cell over time in Fig. \ref{radon1_sat_pres_profiles}. Although we do not have the exact numbers for comparison, a visual inspection indicates that our simulation results match well with those in \cite{BourgeatBenchmark13,Marchand14,Neumann13}.
%\begin{table}[H]
%\centering
%\captionsetup{justification=centering}
%\caption{Performance of Euclid v.s. hypreMGR for mesh size $200 \times 10$} \label{performance_radon1}
%\begin{tabular}{ c | c | c | c | c |}
%\cline{2-5}
%& \multicolumn{2}{|c|}{Linear} & \multicolumn{2}{|c|}{Van Genuchten} \\
%\cline{2-5}
%& ILU(5) & hypreMGR & ILU(5) & hypreMGR \\
%\hline
%Total iters & 10563 & 5432 & & \\
%Total time (s) & 725 & 971 & & \\
%Iters\/Time step & 31.2 & 16.0 & & \\
%\hline
%\end{tabular}
%\end{table}
\par Regarding the setup of hypreMGR in this case, we use three levels of reduction with the restriction and interpolation operators in \cref{non_injective_rp}. For the first level which we need to eliminate the constraints with non-zero diagonal values, a single Jacobi iteration used for the F-relaxation. For the subsequent levels, we apply a singe AMG V(1,1)-cycle for the $A_{ff}$ solve. The coarse grid correction is also solved with one AMG V(2,2)-cycle. \textcolor{black}{With regard to ILU(k), we experiment with different levels of fill and find that ILU(0), ILU(1), etc. would fail to converge for some time step, and ILU(5) is needed for convergence throughout the simulation.} For the nonlinear Van Genuchten model, the new approach requires about \textcolor{black}{34\%} fewer number of iterations and about 17\% less time than ILU(5) for mesh size of $200 \times 10$, as shown in Table \ref{performance_radon1}. The advantage of this approach is much clearer as the problem gets larger (see Table \ref{performance_radon1_400x20}). In this case, hypreMGR takes fewer than half the number of iterations of ILU(5) and requires 40\% less time. Even though the average number of iterations does grow in the case of hypreMGR, it is much less than the rate of ILU(5). Also, we suspect that the mesh may not be large enough for the method to show mesh independence. Similarly, the new approach outperforms ILU(5) in both number of iterations and run time for the linear model.
\begin{table}[H]
\centering
\captionsetup{justification=centering}
\caption{Performance of Euclid ILU v.s. hypreMGR for mesh size $200 \times 10$} \label{performance_radon1}
{\renewcommand{\arraystretch}{1.2}
\begin{tabular}{| c | c | c | c | c | c | c |}
\hline
\multirow{2}{*}{Methods} & \multicolumn{3}{c}{Linear} & \multicolumn{3}{|c|}{Van Genuchten} \\
\cline{2-7}
& Time (s) & LS & LS/NS & Time (s) & LS & LS/NS \\
\hline
ILU(5) & 505 & 22595 & 28.9 & 522 & 24474 & 34.7 \\
hypreMGR & 485 & 16554 & 21.2 & 433 & 15234 & 21.6\\
\hline
\end{tabular}}
\end{table}
\begin{table}[H]
\centering
\captionsetup{justification=centering}
\caption{Performance of Euclid ILU v.s. hypreMGR for mesh size $400 \times 20$} \label{performance_radon1_400x20}
{\renewcommand{\arraystretch}{1.2}
\begin{tabular}{| c | c | c | c | c | c | c |}
\hline
\multirow{2}{*}{Methods} & \multicolumn{3}{c}{Linear} & \multicolumn{3}{|c|}{Van Genuchten} \\
\cline{2-7}
& Time (s) & LS & LS/NS & Time (s) & LS & LS/NS \\
\hline
ILU(5) & 3442 & 42835 & 60.3 & 3467 & 43949 & 64.6 \\
hypreMGR & 2291 & 20408 & 28.7 & 2096 & 19122 & 28.3 \\
\hline
\end{tabular}}
\end{table}

%For all these problems, ILU(5) fails to converge within the allowed maximum number of iterations of 400 for GMRES. We run 50 time steps of $dT = 312.5$ years.
%\begin{table}[H]
%\centering
%\captionsetup{justification=centering}
%\caption{Performance of MGR for different mesh sizes} \label{performance_scaling_radon1}
%{\renewcommand{\arraystretch}{1.1}
%\begin{tabular}{| c | c | c |}
%\hline
%Mesh size & Time & Iters/Time step \\
%\hline
%$2000 \times 10$ & 38.7 & 20.0\\
%$2000 \times 20$ & 56.3 & 20.4\\
%$2000 \times 40$ & 97.4 & 17.2\\
%$2000 \times 80$ & 234.24 & 18.4\\
%$2000 \times 160$ & 337.5 & 23.2\\
%\hline
%\end{tabular}}
%\end{table}
%Note that the number of unknowns is three times the number of mesh points, so the largest problem here has nearly 1 mil unknowns.

\subsection{Three-dimensional Case with Phase Transition}
The domain is a box of dimensions $100\text{m} \times 100\text{m} \times 100\text{m}$. We use a homogeneous permeability field of $\bm{K} = 10^{-14} \text{   m}^2$, which is typical for fresh sandstone (see \cite{Bear72}) that is prevalent in reservoir simulation. The domain is saturated with water, and pure hydrogen is injected into the domain through the boundary of a corner at the bottom. The outlet is set at the opposite corner. The injection rate is $3 \times 5.57\text{ kg/m}^2\text{/year}$ . We run 1 time step of size $dT = 1.825$ days. For the relative permeabilities and capillary pressure models, we use the Van Genuchten model with the same parameters as the example presented in section 6.3. The results are shown in Table \ref{performance_3d}.
%\begin{table}[H]
%\caption{Parameter Values} \label{param_list_3d}
%\centering
%{\renewcommand{\arraystretch}{1.2}
%\begin{tabularx}{0.3\textwidth}{ l X }
%  \hline
%  $P_r$ & $2 \times 10^{6}$ Pa \\
%  $n$ & 1.49 \\
%  $S_{lr}$ & 0.4 \\
%  $S_{gr}$ & 0 \\
%  \hline
%\end{tabularx}}
%\end{table}
\begin{table}[H]
\centering
\captionsetup{justification=centering}
\caption{Performance of hypreMGR for different mesh sizes} \label{performance_3d}
{\renewcommand{\arraystretch}{1.2}
\begin{tabular}{| c | c | c | c | c | c | c |}
\hline
\multirow{2}{*}{Mesh size} & \multicolumn{3}{c}{ILU(5)} & \multicolumn{3}{|c|}{hypreMGR} \\
\cline{2-7} 
& Time (s) & LS & LS/NS & Time (s) & LS & LS/NS \\
\hline
%$10 \times 10 \times 10$ & 12.1 & 308 & 13.4 & 9.0 & 196 & 8.5 \\
$20^3$ & 194.0 & 690 & 25.6 & 121.1 & 267 & 9.9 \\
$40^3$ & 2715.1 & 1470 & 49.0 & 1381.2 & 397 & 13.2 \\
\hline
\end{tabular}}
\end{table}
Since this is a case with phase appearance/disappearance, we use the same setup as in the second example (section 6.2) for hypreMGR. From Table \ref{performance_3d}, the new approach is about 40\% faster in terms of run time, and takes fewer than half the number of iterations of ILU(5) for the mesh size of $20^3$. For the larger mesh size of $40^3$, it is twice faster in terms of run time, and it takes four times fewer the number of iterations of ILU(5). The result indicates that the advantage of hypreMGR clearer as the problem gets larger. Although similar to the result in the previous example, there is an increase in the number of iterations for hypreMGR in the case of the larger mesh, but again the problem is not large enough for us to see the mesh independence result. In fact, we show that this is exactly the case in the next section. And even though the results are not presented here, we note that hypreMGR also outperforms ILU(5) both in run time and number of iterations for the linear model of capillary pressure.

\subsection{Scaling Results}
For the scalability study, we use the same problem setup as in the three-dimensional example in section 6.3. The only difference is in the mesh size. For a strong scaling study, we fix the mesh at $80^3$ (about 1.5 million unknowns) and run the simulation on 8 to 128 cores, each time doubling the number of cores. For weak scaling, we start with a mesh of $40^3$ and then refine the mesh in all directions up to $320^3$, so the largest problem has about 100 million unknowns. We run the problem with 2, 16, 128, and 1024 cores, respectively. The initial time step is $dT = 0.125$ day, and the final time of the simulation is 10 days, except for the case of the largest mesh, which we stop the simulation at 3 days, as we reach the memory limit of the machine.
\par
\begin{figure}[H]
\centering
\includegraphics[width=0.45\textwidth]{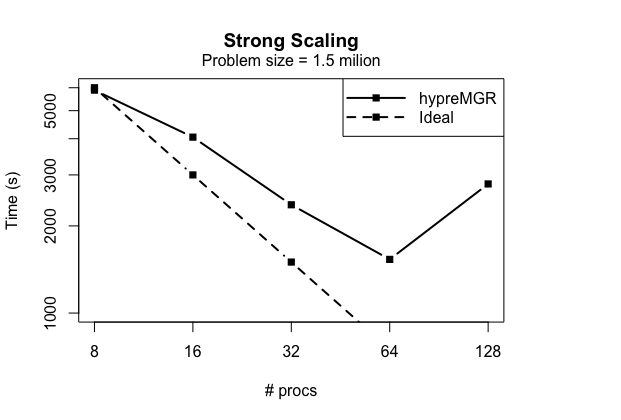}
\includegraphics[width=0.45\textwidth]{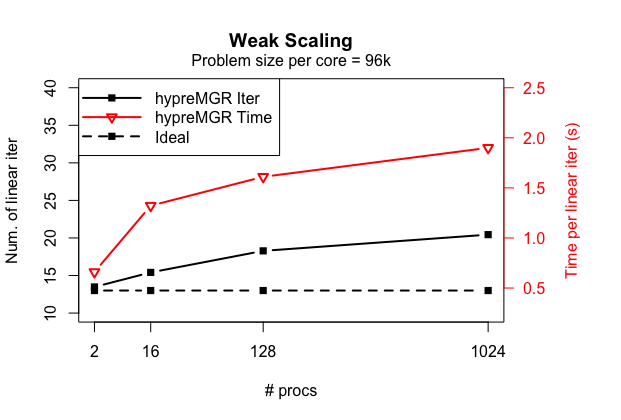}
\caption{Scaling results for MGR.} \label{scaling_results}
\end{figure}

\begin{figure}[H]
\centering
\includegraphics[width=0.45\textwidth]{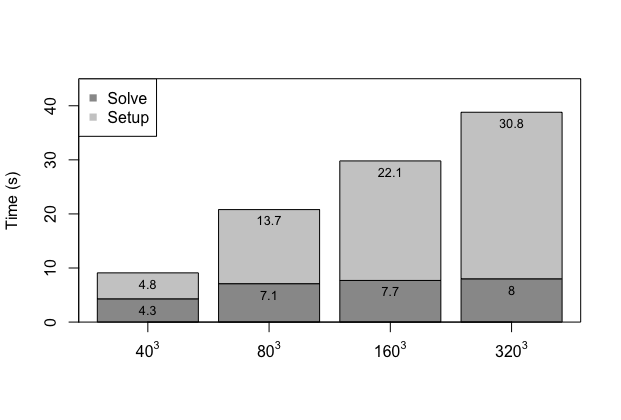}
\caption{CPU time breakdown for linear solve.} \label{ls_time}
\end{figure}

 The results in Fig. \ref{scaling_results} shows that hypreMGR achieves promising results, scaling well up to 64 cores, although it is not quite optimal. From 64 to 128 cores, however, there run time actually increases. This is due to the problem size on each processor getting small (about 12,000 unknowns for 128 cores), and as a consequence, the computation to communication ratio decreases, and that makes the method less efficient. For weak scaling, the performance of hypreMGR is independent of the mesh size. As the problem size gets larger (8 times for each refinement level), the number of linear iterations per nonlinear iterations does not grow significantly, \textcolor{black}{about 14, 19, and 12 percent for 16, 128, and 1024 cores, respectively. The average number of linear iterations also seems to approach a limit, which demonstrates optimal multigrid performance.} Regarding run time, we measure both the setup phase and the solve phase of the algorithm. Since the setup phase requires expensive matrix-matrix multiplications, the total time needed to solve a linear system grows a little faster than the number of iterations. Yet, hypreMGR still achieves near optimal scalability.
\textcolor{black}{Fig. \ref{ls_time} focuses on the time of the linear solve, splitting into the setup and solution phases. It is clear that solve phase achieves optimal scalability as the time needed to iterate to convergence stays nearly constant for mesh sizes $80^3$, $160^3$, and $320^3$. In contrast, the setup phase, which includes constructing $R$ and $P$, computing the coarse grids using the matrix-matrix product $RAP$, and all the AMG setup for the coarse grid as well as the F-relaxation, does not scale very well. This is likely an implementation problem, and it can be improved in the future.}

In terms of memory storage, like AMG methods, the MGR approach requires the storage of the restriction, interpolation, and coarse grid operators at every level. However, in addition to these operators, the MGR approach also needs to store the $A_{ff}$ matrices at the levels which scalar AMG is used for the F-relaxation step. When aggressive coarsening is performed in standard AMG methods, the size of the coarse grid can be significantly reduced after the first level. With MGR, the size of the matrix is reduced by at most a third, since the reduction is dictated by the block structure of the system, rather than the heuristics used in AMG.

\section{Conclusion}
\label{sec:conclusion}
We have presented a preconditioning strategy for solving the linear systems that arise from the solution of multiphase multicomponent porous media flow with phase transitions. To account for the phase transitions, the problem is formulated as a nonlinear complementarity problem, and solved using the semi-smooth Newton technique. 

The proposed preconditioner is based on the multigrid reduction technique, which generalizes traditional two-stage preconditioners in a natural multigrid framework. In this work, we extend a previously developed two-grid strategy to a multilevel reduction strategy that accounts for the transitions in the phases of the primary variables. We have demonstrated the performance of the preconditioner on classic benchmark problems presented in the literature, and show the parallel efficiency of the linear solver on large-scale problems. The numerical results indicate optimal scalability and robust performance of the MGR preconditioner, which is important for real-field simulations.

We observed that depending on the properties of the capillary pressure model used, a different solver could be used for the F-relaxation phase of the preconditioner. When the model is convection-dominated, a simple relaxation scheme is sufficient for F-relaxation. However, when the model is diffusion-dominated, relaxation alone is not sufficient, and a more robust solver is required for F-relaxation. In our experiments, we used AMG for such problems. However, this may be excessive and in some cases, inefficient. \textcolor{black}{For applications with phase transitions, the fronts along which the transitions occur can be small compared to the entire domain. As a result, using the same strategy for F-relaxation at these intermediate solves can be inefficient since communication dominates computation at this point. Allowing different strategies to be employed, as dictated by the physics, can be a more efficient strategy. For example, using a single-level relaxation strategy instead of a multilevel (v-cycle) technique could be more appropriate. We are exploring this idea, in addition to aggressive coarsening strategies to improve parallel efficiency.}

\textcolor{black}{Future applications of interest for the MGR solver include applications with multiple phases, poromechanics, and applications with fractures and thermal properties. The MGR framework is general enough to handle these applications as a ``black-box" solver, and can also serve as a basis for building good physics-based preconditioners. However, more work may be required to improve solver performance for these complex applications. We are exploring new strategies for building interpolation and restriction operators so that the final coarse grid system is a good approximation to a pressure system (or has elliptic\/ M-matrix properties) that is amenable to solution by AMG. We are also considering incorporating structure information within a semi-structured framework to develop a robust solver that can effectively handle grid anisotropy for complex geometries.
}

\bibliographystyle{plain}
\bibliography{copper17}

\end{document}

%% file: ex_shared.tex
\usepackage{lipsum}
\usepackage{amsfonts}
\usepackage{graphicx}
\usepackage{epstopdf}
\usepackage{algorithmic}
\ifpdf
  \DeclareGraphicsExtensions{.eps,.pdf,.png,.jpg}
\else
  \DeclareGraphicsExtensions{.eps}
\fi

\renewcommand{\vec}[1]{\mathbf{#1}}

\newcommand{\TheTitle}{Algebraic Multigrid Preconditioners for Multiphase Flow \\ in Porous Media with Phase Transitions} 
\newcommand{\TheAuthors}{Quan M. Bui, Lu Wang, and Daniel Osei-Kuffuor}

\newcommand{\HeaderTitle}{Algebraic Multigrid Preconditioners for Multiphase Flow}

\headers{\HeaderTitle}{\TheAuthors}

\title{{\TheTitle}\thanks{This work was performed under the auspices of the U.S. Department of Energy by Lawrence Livermore National Laboratory under Contract DE-AC52-07NA27344.}}

\author{
  Quan M. Bui\thanks{Corresponding author. Applied Math, Stats, and Scientific Computation, University of Maryland, College Park, MD {\texttt email:} \email{mquanbui@math.umd.edu}.}
  \and
  Lu Wang \thanks{Center for Applied Scientific Computing, Lawrence Livermore National Laboratory, Livermore, CA. {\texttt email:} \{\email{wang84, oseikuffuor1}\}\email{@llnl.gov}.}
  \and
  Daniel Osei-Kuffuor \footnotemark[3]
}

\usepackage{amsopn}